\documentclass{article}
\usepackage{a4wide}

\usepackage{lineno}
\usepackage{microtype}
\usepackage{subcaption}
\usepackage[onehalfspacing]{setspace}
\usepackage{epsfig}

\usepackage{tikz}
\usetikzlibrary{patterns}
\usepackage{pgfplots}
\pgfplotsset{compat=newest}
\usepackage{comment}
\usepackage{amssymb}
\usepackage{amsmath}
\usepackage[noabbrev,capitalise,english]{cleveref} 

\graphicspath{{./figs/}}

\usepackage{color}

\usepackage{bm}

\newcommand{\revision}[1]{\textcolor{red}{#1}}
\newcommand{\franz}[1]{\textcolor{blue}{#1}}
\newcommand{\phuoc}[1]{\textcolor{red}{#1}}

\newcommand{\michel}[1]{\textcolor{red}{#1}}

\usepackage{algorithm, algorithmicx, algpseudocode}

\usetikzlibrary{shapes}
\usetikzlibrary{plotmarks}
\usepackage{float}

\newcommand{\be}{\begin{equation}}
\newcommand{\ee}{\end{equation}}
\newcommand{\beq}{\begin{eqnarray}}
\newcommand{\eeq}{\end{eqnarray}}
\newcommand{\beqs}{\begin{eqnarray*}}
\newcommand{\eeqs}{\end{eqnarray*}}
\newcommand{\et}{\end{theorem}}
\newcommand{\bex}{\begin{example}}
\newcommand{\eex}{\end{example}}
\newcommand{\br}{\begin{remark}}
\newcommand{\er}{\end{remark}}
\newcommand{\bc}{\begin{corollary}}
\newcommand{\ec}{\end{corollary}}
\newcommand{\bl}{\begin{lemma}}
\newcommand{\el}{\end{lemma}}
\newcommand{\bp}{\begin{proposition}}
\newcommand{\ep}{\end{proposition}}
\newcommand{\bd}{\begin{definition}}
\newcommand{\ed}{\end{definition}}
\newcommand{\bas}{\begin{assumption}}
\newcommand{\eas}{\end{assumption}}

\newcommand{\R}{\mathbb{R}}

\newcommand{\PP}{\mathbb{P}}

\def\u{{\bf u}}

\def\n{{\bf n}}

\def\B{{\bf B}}
\def\F{{\bf F}}
\def\0{{\bf 0}}

\newcommand{\bPi}{\bm{\Pi}}
\newcommand{\bN}{\bm{N}}

\newcommand{\bB}{\bm{B}}
\newcommand{\bT}{\bm{T}}
\newcommand{\bzero}{\bm{0}}
\newcommand{\bu}{\bm{u}}
\newcommand{\bpp}{\bm{p}}

\newcommand{\bE}{\bm{E}}

\newcommand{\bC}{\bm{C}}
\newcommand{\bI}{\bm{I}}
\newcommand{\bF}{\bm{F}}
\newcommand{\bV}{\bm{V}}
\newcommand{\bQ}{\bm{Q}}
\newcommand{\bq}{\bm{q}}
\newcommand{\bv}{\bm{v}}
\newcommand{\bw}{\bm{w}}
\newcommand{\bz}{\bm{z}}
\newcommand{\bR}{\bm{R}}
\newcommand{\bJ}{\bm{J}}

\newcommand{\eg}{e.g.} 

\usepackage[english]{babel}

\begin{document}

 \title{Enhancing Biomechanical Simulations Based on A Posteriori Error Estimates: The Potential of Dual Weighted Residual-Driven Adaptive Mesh Refinement} 

 \author{Huu Phuoc Bui\footnote{Ansys , Lyon, France, \texttt{hphuoc.bui@gmail.com}}, 
 Michel Duprez\footnote{MIMESIS, MLMS, Inria Nancy Grand-Est, Université de Strasbourg, Strasbourg, France, \texttt{michel.duprez@inria.fr}},
Pierre-Yves Rohan\footnote{Institut de Biomécanique Humaine Georges Charpak, Arts et Métiers Institute of technology, Paris, France, \texttt{pierre-yves.rohan@ensam.eu}},
Arnaud Lejeune\footnote{Department of Applied Mechanics, FEMTO-ST Institute, University of Franche-Comte, UMR 6174 CNRS, Besan\c{c}on, France},\\
Stéphane P.A. Bordas\footnote{Institute of Computational Engineering, Department of Engineering, Luxembourg, \texttt{stephane.bordas@gmail.com}},
Marek Bucki\footnote{TwInsight, Grenoble, France, \texttt{marek.bucki.pro@gmail.com}},
Franz Chouly\footnote{Institut de Mathématiques de Bourgogne, Université de Bourgogne, Dijon, France}
\footnote{Mathematical Modeling and Department of Mathematical Engineering, University of Chile, Santiago, Chile}
\footnote{Departamento de Ingeniería Matemática, CI2MA, Universidad de Concepción, Concepci\`on, Chile, \texttt{franz.chouly@gmail.com}}
}

\maketitle

\abstract{
The Finite Element Method (FEM) is a well-established procedure for computing approximate solutions to
deterministic engineering problems described by partial differential equations. FEM produces discrete approximations of the solution with a
discretisation error that can be an be quantified with \emph{a posteriori} error estimates. 
The practical relevance of error 
 estimates for biomechanics problems, especially for soft tissue where the response is governed by large strains, is rarely addressed. 
In this contribution, we propose an implementation of \emph{a posteriori} error estimates targeting a user-defined quantity of interest, using the Dual Weighted Residual (DWR) technique tailored to biomechanics. The proposed method considers a general setting that encompasses three-dimensional geometries and model non-linearities, which appear in hyperelastic soft tissues. We take advantage of the automatic differentiation capabilities embedded in modern finite element software, which allows the error estimates to be computed generically for a large class of models and constitutive laws. First we validate our methodology using experimental measurements from silicone samples, and then illustrate its applicability for patient-specific computations of pressure ulcers on a human heel.
}


\section{Introduction}
\label{sec:Introduction}

The importance of finite element analyses (FEA) for biomechanical investigations has increased considerably worldwide in recent years. Such finite element models are widely employed to investigate both the underlying mechanisms that drive normal physiology of biological soft tissues and the mechanical factors that contribute to the onset and development of diseases such as tumour growth \cite{Yankeelov187ps9,urcun_2023}),  atherosclerosis or aneurysms \cite{romo_vitro_2014}, or multilevel lumbar disc degenerative diseases \cite{schmidt_effect_2012}, to name a few. Finite element models are also valuable tools that contribute to the development of medical devices such as, for example, vascular stent-grafts \cite{perrin_patient-specific_2015}, and have the potential to improve prevention strategies \cite{Kranke2007,Trabelsi2015patient,luboz_personalized_2017}, surgical planning \cite{Vannier1984,Selle2002,buchaillard_simulations_2007}, pedagogical simulators for medical training \cite{Kuhnapfel2000,courtecuisse2014real} and guidance of surgeons during interventions \cite{Warfield2002,Carter2005}. A survey of applications using simulation modelling for healthcare sector can be found for instance in \cite{Mielczarek2012}.

In this context{,} one major issue is meshing, since the reliability of the predicted {mechanical response} arising from computer simulation {heavily relies} on the quality of the underlying finite element mesh
{\cite{gratsch_posteriori_2005}}.
The patient-specific mesh has to be built from segmented medical images (CT, MRI, ultra-sound){,} and {has} to conform to anatomical details with potentially complex topologies and geometries {\cite{bijar2016atlas}}. 
In general the quality of a given mesh is assessed through purely geometrical criteria, that allow in some way to quantify the distortion of the geometry of the elements 
\cite{bucki-2011}.
Beyond mesh quality, mesh density is another, related, parameter which must be controlled during biomechanics simulations. 
Moreover solutions must be obtained on commodity hardware within clinical time scales: milliseconds (for surgical training);  minutes (for surgical assistance); hours (for surgical planning).
As a result one question that always arises in practice is: ``Given a tolerable error level, what is the coarsest possible mesh which will provide the required accuracy?''
This leads to the notion of ``mesh optimality'', which is achieved for an optimal balance between the accuracy in a given quantity of interest to the user and the associated computational cost.

In this paper, we investigate the capability of
\emph{a posteriori} error estimates \cite{ainsworth-oden-2000,verfurth-2013} to provide useful information about the {discretization error}, \emph{i.e.},
the difference between the finite element solution and the exact solution of the same boundary value problem on the same geometry.
\emph{A posteriori} error estimates are quantities computed 
from the  numerical solution, that indicate the magnitude of the local error. These estimates are at the core of mesh adaptive techniques \cite{nochetto-2009}. 
Many {\it a posteriori} error estimation methods have been developed in the numerical analysis community. These methods have different theoretical and practical properties. However, despite their great potential, 
error estimates have rarely been considered for patient-specific finite element simulations {in the biomechanical community}. 

To the best of our knowledge, the first works that address this issue are 
\cite{Bui2016_TBME,Bui2017_error_brain_ijnmbe}, which study the discretization error (based on energy norm) of real-time simulations using the recovery-based 
technique of Zienkiewicz and Zhu \cite{zienkiewicz1992superconvergent}. This approach is inexpensive and allows to deal with real-time simulations. However, the error in energy norm might not provide useful information for applications where one is interested in the error of a real physical quantity of interest. To overcome this difficulty, estimates based on duality arguments are common for {\it a posteriori} error estimation, see \eg{} \cite{becker-rannacher-2001,Heuveline2003,paraschivoiu-peraire-patera-1997,prudhomme-oden-1999,maday-patera-2000,giles-suli-2002}.
A preliminary study has been carried out previously in this direction in \cite{duprez2020} by the authors of the present paper. This study makes use of the Dual Weighted Residuals (DWR) method, as presented in \cite{becker-rannacher-2001}.
Let us recall that the main idea of this method is to solve a dual problem, the solution of which is used as a weight that indicates locally the sensitivity of the quantity of interest for each cell-wise contribution to the discretization error. However, this aforementioned study is limited to a simplified setting, since it was aimed at giving preliminary insights and at addressing the first technical difficulties. The modelling of soft tissues in \cite{duprez2020} is indeed restricted to two-dimensional linear elasticity (plane strain) {problems} and to a quantity of interest that should depend linearly of the displacement.

As a result, the main goal of this paper is to handle a setting much closer to current practice in soft tissue simulation. For this purpose, we consider three-dimensional (passive) hyperelastic soft tissue 
and arbitrary quantities of interest, that may depend non-linearly on the displacement. The DWR method is very well adapted to this setting, as it was designed originally for non-linear problems \cite{becker-rannacher-2001} (see also \cite{larsson-2002,whiteley-tavener-2014,gonzalezestrada2014,granzow2018} 
for applications in non-linear elasticity). Nevertheless, one major issue for its application is the practical calculation of the dual solution, which involves the derivatives of the primal weak form and of the quantity of interest (this problem does not appear in the linear setting). This formal derivation can be intricate for soft tissue models built from complex hyperelastic constitutive laws. To handle this issue, we take advantage of the capabilities of modern finite element softwares such as FEniCS or GetFEM++, that integrate automatic symbolic differentiation. Not only it makes easier the implementation of the DWR method for error estimation, but also it requires no real effort if the constitutive law is changed.
Then, we validate the methodology using experimental data obtained from in-vitro study of silicone samples \cite{meunier2008mechanical} and show its potential interest on 
an example coming from patient-specific simulation. 

%

This paper is organized as follow{s}. In Section~\ref{sec:Methods}, we describe the hyperelastic setting for 
passive soft tissue, the corresponding finite element discretization, the DWR \textit{a posteriori} error estimation as well as the algorithm for mesh refinement. In Section~\ref{sec:Results}, we illustrate the methodology for different test-cases. 
The results are discussed in Section~\ref{sec:Conclusions}.


\section{Methods}
\label{sec:Methods}

First we present the model problem, then the finite element discretization and finally the error estimation and mesh refinement techniques.

\subsection{Problem setting: incompressible hyperelastic soft tissue}


We consider an (incompressible) hyperelastic body in a reference configuration denoted by $\Omega$, an open subset of $\R^3$, and subjected to a given body force $\bB$. The unknown displacement field and the unknown static pressure are denoted by $\bu$ and $\bpp$, respectively.
The deformation gradient is denoted $\bF$, with $\bF := \bI + \nabla_X \bu$, where $\bI$ stands for the identity matrix, and $\nabla_X$ denotes the gradient with respect to coordinates in the reference configuration. 
The first Piola-Kirchhoff stress tensor denoted $\bPi$ is derived from the hyperelastic strain-energy density function $W$, which depends on the  displacement field $\bu$ and the pressure $\bpp$, as follows:
\begin{equation}
\bPi = \frac{\partial W}{\partial \bF}.
\end{equation}
In the present paper, we will consider different incompressible material models corresponding to some strain-energy densities, namely:
\begin{enumerate}
\item \textbf{Mooney-Rivlin model (see  \cite{Mooney1940}):}
\begin{equation}\label{eq:mooney}
W := c_{10} (J_1 - 3) + c_{01} (J_2 - 3) -\bpp(\mbox{det}(\bC)-1) 
\end{equation}
\item \textbf{Gent model (see  \cite{gent2}) :}
\begin{equation}\label{eq:gent}
W := \frac{-EJ_m}{6}\ln\left(1-\frac{J_1-3}{J_m}\right) -\bpp(\mbox{det}(\bC)-1) 
\end{equation}
\item \textbf{Haines-Wilson model (see  \cite{james1975strain}):}
\begin{multline}\label{eq:HW}
W := c_{10} (J_1 - 3) + c_{01} (J_2 - 3) +c_{20} (J_1 - 3)^2 + c_{02} (J_2 - 3)^2\\+ c_{30} (J_1 - 3)^3+c_{11} (J_1 - 3)(J_2 - 3) -\bpp(\mbox{det}(\bC)-1) 
\end{multline}
\end{enumerate}
where  $\bC (:= \bF^T \cdot \bF)$ denote the right Cauchy-Green tensor,
 $J := \textrm{det } \bF$ the jacobian of the deformation, where $I_1 := \textrm{trace } \bC$, $I_2 := \frac12 ( (\textrm{tr } \bC)^2 - \textrm{tr } (\bC \cdot \bC))$, $J_1 := I_1 J^{-\frac23}$ and $J_2 := I_2 J^{-\frac43}$ are invariants associated to the deformation and 
$c_{ij}$, $J_m$ and $E$ are some coefficients which will be given in Table \ref{tab:parameter_law}. 
For the sake of simplicity, the boundary $\partial \Omega$ of $\Omega$ is partitioned into two subsets $\Gamma_D$ and $\Gamma_N$, and we apply a prescribed displacement $\bu = \bu_D$ on $\Gamma_D $ and a given  force $\bT$ on $\Gamma_N$. 

Let us introduce the virtual works associated to the internal and external forces:
\begin{equation*} 
A(\bu,\bpp; \bv,\bq) := \int_{\Omega} \bPi(\bu,\bpp): \nabla_X \bv \, \mathrm{d} \Omega
+ \int_{\Omega} (1-\mathrm{det}(\bC))\bq \, \mathrm{d} \Omega, \qquad
L(\bv) := \int_{\Omega} \bB \cdot \bv \, \mathrm{d} \Omega   +
\int_{\Gamma_N} \bT \cdot \bv  \,\mathrm{d} \Gamma,
\end{equation*}
where $\bu$ and $\bv$ are admissible displacements and $\bpp$ and $\bq$ are admissible pressures.
The hyperelastic problem in weak form reads
\begin{equation}
\left\{
\begin{array}{l}
\text{Find a displacement }\bu, 
\text{with } 
\bu = \bu_D \text{ on } \Gamma_D
\text{ and a pressure }\bpp 
\text{ such that } \\\noalign{\smallskip}
A(\bu,\bpp; \bv,\bq) = L(\bv), \, \forall (\bv,\bq), 
\, \bv = \mathbf{0} \text{ on } \Gamma_D,
\end{array}
\right.
\label{eq:primal_weak}
\end{equation}
%
Let $\mathcal{K}_h$ be a mesh of the domain $\Omega$.
Let us denote by $\bV_h \times \bQ_h$ the finite element pair that makes use of the lowest-order Taylor-Hood finite elements on $\mathcal{K}_h$ (continuous piecewise polynomials of order $2$ for the displacement and of order $1$ for the pressure). The finite element method to solve our hyperelastic problem reads 
%
\begin{equation}
\left\{
\begin{array}{l}
\text{Find a displacement }\bu_h \in \bV_h \text{, with } \bu_h = \bu_D^h \text{ on } \Gamma_D \text{ and a pressure }\bpp_h\in\bQ_h\text{ such that } \\\noalign{\smallskip}
A(\bu_h,\bpp_h; \bv_h,\bq_h) = L(\bv_h), \, \forall (\bv_h,\bq_h) \in \bV_h^0\times\bQ_h,
\end{array}
\right.
\label{eq:FEM_primal_weak}
\end{equation}
where $\bV_h^0$ is composed by the functions of $\bV_h$ vanishynig on $\Gamma_D$ and where $u_D^h$ is a finite element approximation of $u_D$, obtained for instance by Lagrange interpolation or by projection.

Note that the above choice of Taylor-Hood finite elements is for the sake of simplicity, and that the methodology described below for mesh refinement can be extended rather straightforwardly to other conforming variational discretization techniques, provided they ensure a stable and accurate approximation of the finite elasticity equations, and that they allow to split the residual as a sum of local contributions. For instance, any other infsup stable pair of finite elements (mini-element, P2-iso-P1) on simplicial, tensor-product or mixed meshes can be considered.

\subsection{Quantity of interest}
The objective of goal-oriented error estimation is to estimate the error of the finite element solution 
for a user-defined quantity, possibly different from the  energy norm, see \cite{becker-rannacher-1996,becker-rannacher-2001,giles-suli-2002,gonzalezestrada2014,ODEN2001_Goal}. The idea of this technique comes from the fact that one would like to analyse the error of a predefined target quantity since in certain circumstances the global error norm may not be useful. 

Let $Q(\bu,\bpp)$ be a quantity of biomechanical interest computed from the exact solution $(\bu,\bpp)$ of the problem \eqref{eq:primal_weak}, 
with smooth, but possibly nonlinear, dependence on $(\bu,\bpp)$.
%
The aim is to estimate the error for the quantity of interest $$|Q(\bu,\bpp)-Q(\bu_h,\bpp_h)|$$ where $(\bu_h,\bpp_h)$ is the approximated finite element solution of \eqref{eq:FEM_primal_weak}. 

\subsection{Dual problem for computing the weights}

One of the main ingredients of the DWR method is to solve an adjoint problem to extract information from the quantity of interest, and get feedback about the regions where it is more, or less, influenced, by the approximation error. As a result, we compute a pair of dual variables
$(\bz_h,\bw_h)$ by solving 
\begin{equation}
\left\{
\begin{array}{l}
\text{Find } ({\bz}_h,\bw_h) \in {\bV}_h^0  \times \bQ_h \text{ such that } \\\noalign{\smallskip}
(A')^*(\bu_h,\bpp_h| {\bz}_h,  {\bw}_h;{\bv}_h, {\bq}_h) = Q'(\bu_h,\bpp_h|{\bv}_h,\bq_h) \quad \forall ({\bv}_h,\bq_h) \in {\bV}_h^0  \times \bQ_h,
\end{array}
\right.
\label{eq:discrete_dual_problem}
\end{equation}
where $A'$ and $Q'$ denote the Fr\'echet derivative of $A$ and $Q$, respectively,
and $(A')^*$ is the adjoint form of $A'$. 
Then we extrapolate the solution in a finite element space of higher polynomial degree following \cite[Fig. 5.1]{rognes-logg-2013}. This function weights the residual in our estimator.

Remark that the dual problem \eqref{eq:discrete_dual_problem} is linear, so solving it is not expensive in comparison to \eqref{eq:primal_weak}.
For model problems such as \eqref{eq:primal_weak}, and some expressions of $Q$, the practical calculation of $A'$ and $Q'$ can be tough. For this purpose, we take advantage of the capabilities of automatic symbolic differentiation embedded into modern finite element software such as FEniCS or GetFEM++. 
Furthermore, this feature makes possible some genericity in the implementation: virtually nothing has to be changed in the program if the hyperelastic constitutive law is modified.

\subsection{The representation formula of Becker and Rannacher}

We introduce $r(\bu_h,\bpp_h;\bv,\bq)$ the residual of Problem \eqref{eq:FEM_primal_weak} as
\begin{equation}
r(\bu_h,\bpp_h;\bv,\bq) = L(\bv,\bq) - A(\bu_h,\bpp_h;\bv,\bq) \qquad \forall (\bv,\bq) \in \bV\times\bQ.
\end{equation}
This, roughly speaking, quantifies how well the hyperelasticity equations are approximated (it should tendsto zero if the mesh is uniformly refined).
Thanks to the dual system \eqref{eq:discrete_dual_problem}, we obtain expression of the error on $Q$ 
as a best approximation term involving the residual and the (exact) dual solution
(see \cite[Proposition 2.3]{becker-rannacher-2001}): 
\begin{equation}
Q(\bu,\bpp) - Q(\bu_h,\bpp_h) = \min_{(\bv_h,\bq_h) \in \bV_h\times\bQ_h} r(\bu_h,\bpp_h; \bz - \bv_h, \bw - \bq_h) + R_m \label{eq:representation_DWR}
\end{equation}
where 
$R_m$ is the high order remainder related to the error caused by the linearization of the nonlinear problem (the precise expression of which can be found in \cite{becker-rannacher-2001}). 
In practice, this quantity is, hopefully, negligible.
Note at this stage that there are various possibilities to represent the error on $Q$, which are detailed in \cite{becker-rannacher-2001}, and for instance, in \cite{rognes-logg-2013}, the authors make use of another representation formula (Proposition 2.4) which is then approximated. 

Proceeding as usual in \emph{a posteriori} error estimation, \emph{i.e.,}
after performing integration by parts on the residual $r$, we localize the different contributions to the goal-oriented error as follows:
\begin{equation}
\lvert Q(\bu,\bpp) - Q(\bu_h,\bpp_h) \rvert \leq \sum_{K \in \mathcal{K}_h} \eta_K{((\bu_K,\bpp_K),(\bz_K,\bw_K))} + H.O.T.
\end{equation}
In the above expression, $K$ denotes any cell of the mesh $\mathcal{K}_h$, and expressions such as $\bu_K$ denote the local restriction of the finite element variable $\bu_h$ ot the cell $K$. Moreover $H.O.T.$ denotes high order terms, that are not considered in the implementation. In Appendix are provided the detailed expression of $\eta_K$.

\subsection{Adaptive mesh refinement}

Using the error estimate on $Q$, we implement a standard procedure for mesh refinement. 
As described in \cref{algo:adaptivity}, we start with an initial mesh called $mesh_i$, and by providing a guest solution $\bu_i^{(0)}$, the nonlinear primal problem can be solved using the Newton's method (see \cref{algo:newton_method}). Once accepting $\bu_i$ as the solution of the primal problem, solving the discrete dual problem (see \cref{algo:dual_problem}) provides the dual solution $\bz_i$ of the actual mesh $mesh_i$. The estimator $\eta_K$ is then computed providing the primal and dual solutions $\bu_i$ and $\bz_i$, respectively. From the estimator, different strategies can be used to mark the elements whose \emph{error} is high. In this paper, we use the D\"orfler marking strategy \cite{dorfler1996} (see Algorithm \ref{algo:Dorfler_making}).
The detailed algorithms are given in the Appendix.


\section{Results} 
\label{sec:Results}
In this section, we present the performance of the DWR method in controlling the discretisation error in simulations employing hyperelastic models. We will consider two test cases. In the first one, we will compare different constitutive laws and compute the model error thanks to experimental data. We will then highlight the performance of the DWR strategy to reduce the discretisation error. The second test case will be only in silico on a heel geometry.
The simulations have been realized thanks to the \texttt{python} library \texttt{FEniCS} and the code is available online\cite{Duprez2023}.




\subsection{First test case : silicone samples}

\begin{figure}[H]
     \centering
~\hfill  \begin{minipage}{0.25\linewidth}
  \includegraphics[width=1\textwidth]{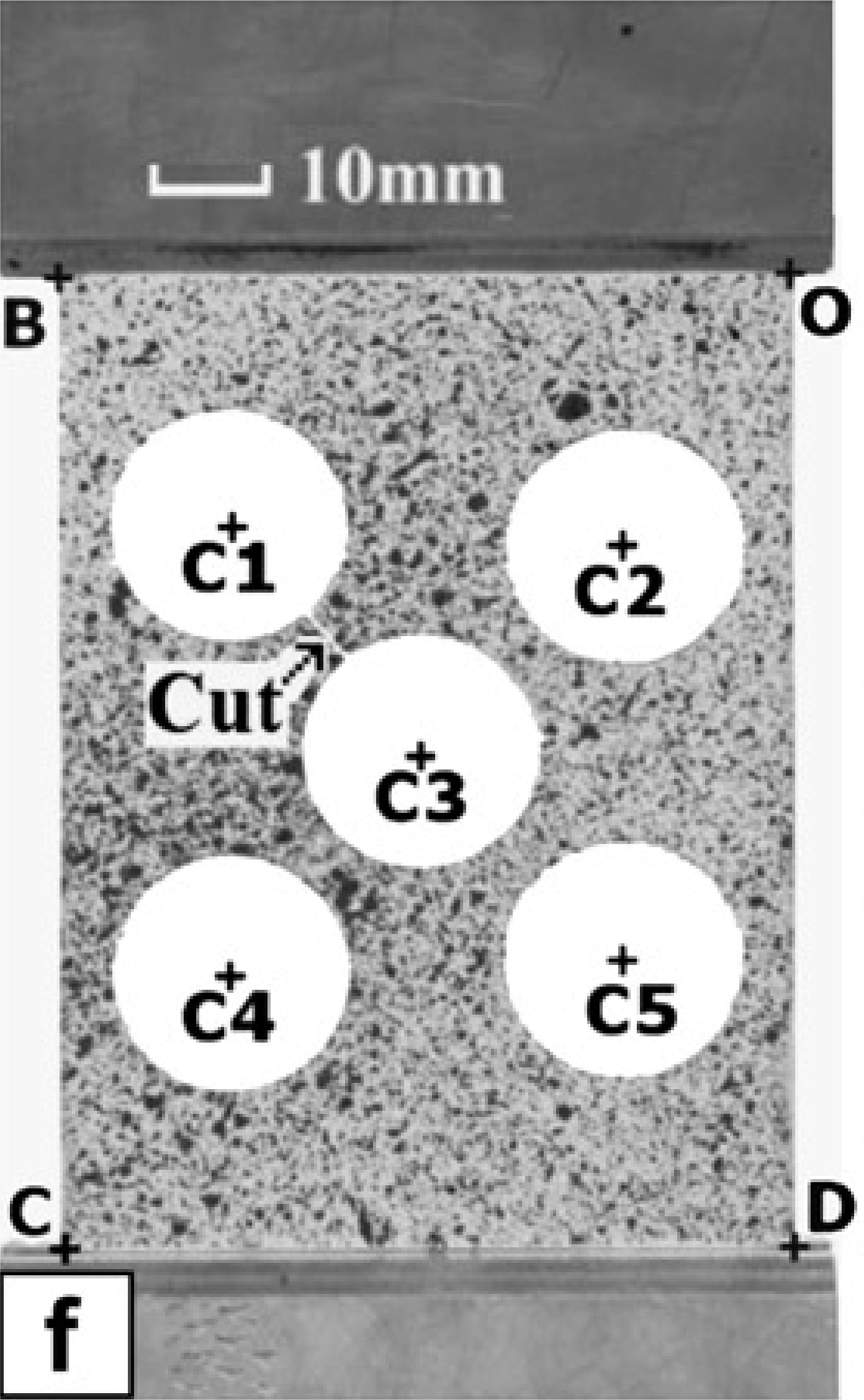}\vspace*{-0cm}
\end{minipage}\hfill
    \begin{minipage}{0.3\linewidth}
     \begin{tabular}{|c|c|c|}\hline
          &X (mm)&Y (mm)  \\\hline
          O& 0&0\\ \hline
        B& -62&0\\ \hline
      C& -62&-82.5\\ \hline
D& 0&-82.5\\ \hline
          $C_1$& -47.5&-21.4\\ \hline
    $C_2$& -14.0&-23.0\\ \hline
      $C_3$& -31.5&-41.0\\ \hline
        $C_4$& -47.7&-59.1\\ \hline
        $C_5$& -14.5&-58.0\\ \hline
     \end{tabular}\end{minipage}\hfill
     \begin{minipage}{0.25\linewidth}
     \includegraphics[width=1\linewidth]{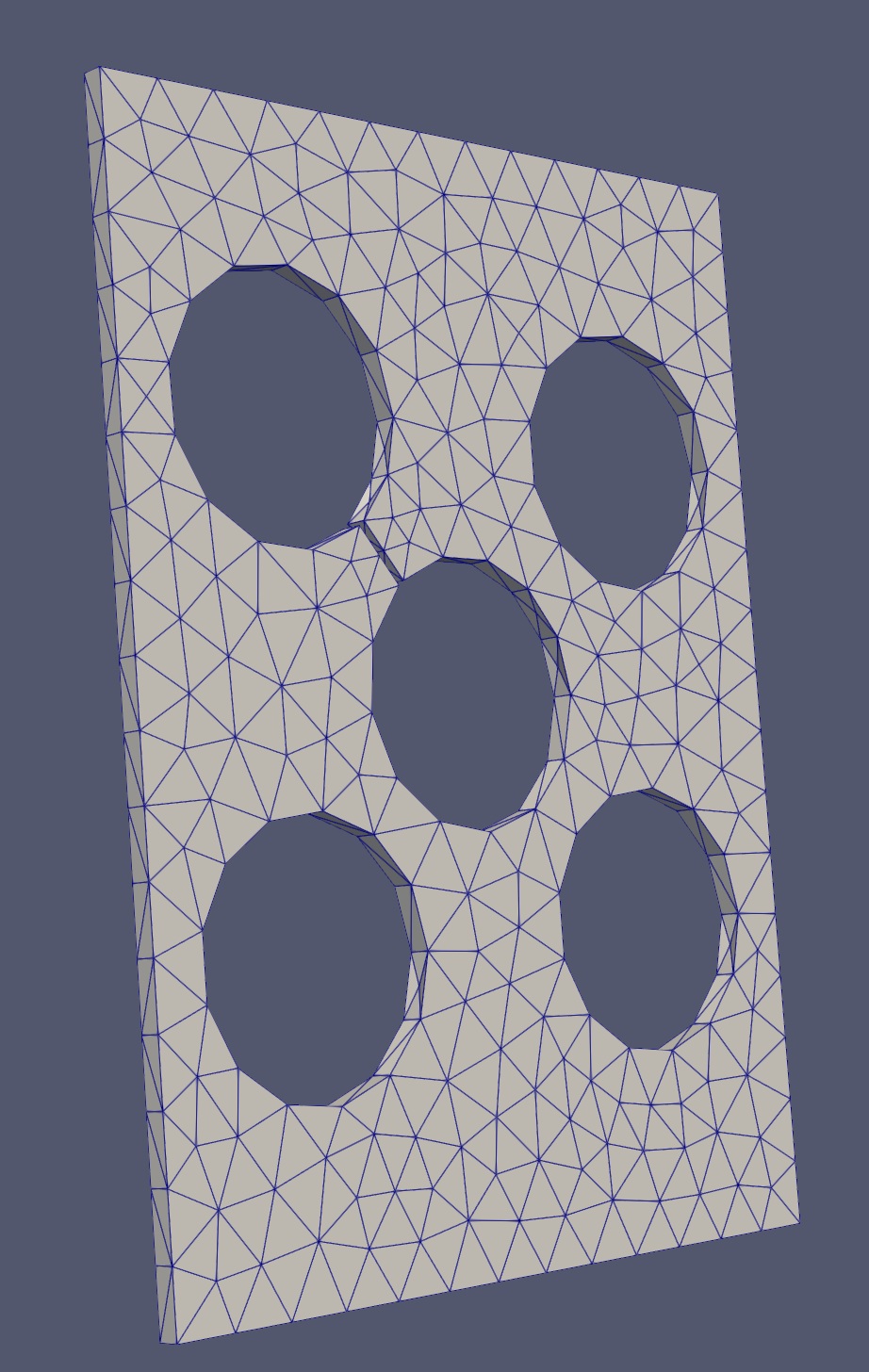}\vspace*{-0.3cm}
     \end{minipage}\hfill~
     \caption{Geometry of the silicone sample, position of the holes and initial mesh.}
     \label{fig:geometry}
    \end{figure}

The experimental procedure is briefly recalled here for the sake of clarity. For more information, the reader is referred to \cite{meunier2008mechanical}.
Simple tensile tests are performed on dumbbell shaped samples of silicone rubber (RTV 141) having an initial gauge length $l_{0}$ of 82.5 mm, a gauge width $b_{0}$ of 61.5 mm and a gauge thickness $e_{0}$ of 1.75 mm.  The sample contains five holes of diameter 20 mm
and the position of the centres of the holes and the corners are given in Figure \ref{fig:geometry}. There is also a cut between the circles C1 et C3.

Tested samples are deformed using a universal mechanical testing machine (MTS 4M) (see \cite{meunier2008mechanical}).
Dirichlet boundary conditions are imposed on the bottom edge of the dumbell silicone sample. On the left and right boundaries, we impose a homogeneous Neumann boundary condition ($\F = \0$). In the experiment, the following Neumann boundary condition is imposed on the top edge: 
$\F = \left (f_A / (b_{0} \times e_{0} )\right ) \n$ 
such that $\int_{\Gamma_N} \F \cdot \n \, \mathrm{ds} = 20$ N.
This force implies an observed vertical displacement of 57.3mm. In the simulations, we do it the other way round: we fix the bottom and we impose a displacement of 57.3mm on the top. We guess the corresponding traction force on the top boundary.
Thus,we will consider the following quantities of interest 
$$Q(\bu,\bpp)=\int_{\mbox{top}}(\bPi(\bu) \cdot \bN)\cdot \bN \mathrm{d}s,$$
where the integral is taken on the top of the silicon band.

\begin{figure}
    \centering
    \includegraphics[width=0.3\linewidth]{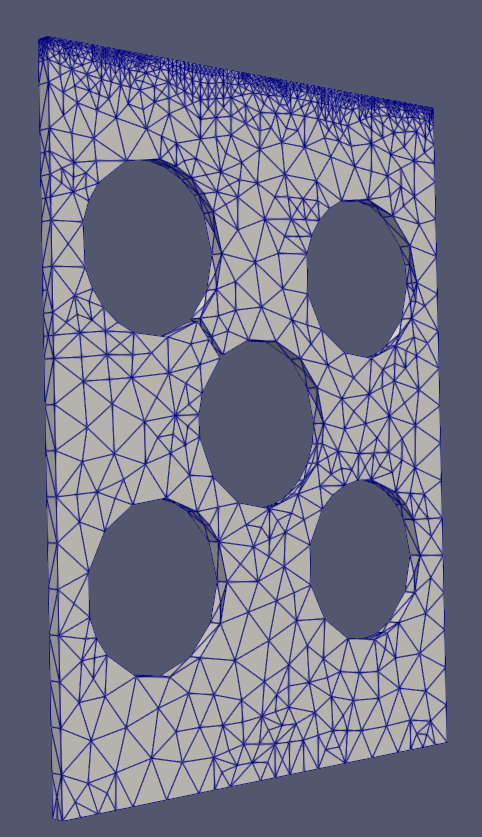}
        \includegraphics[width=0.3\linewidth]{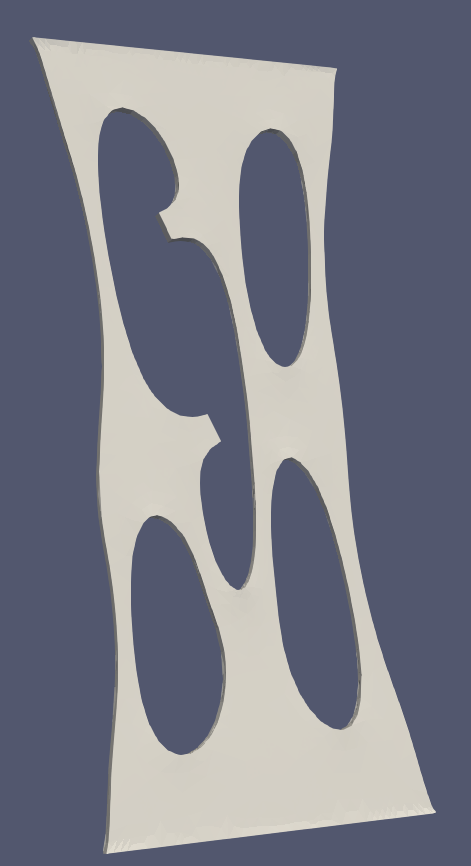}
    \caption{lar}       
        \caption{First test case (silicone sample).  Refined mesh for the Haine-Wilson model (left); Deformed geometry (right)}
    \label{fig:deformed mesh}
\end{figure}

Table \ref{tab:parameter_law} recalls the value of the constitutive parameters used in the simulations.
In Table \ref{tab:law}, we compare the model and discretisation error  for  the constitutive laws of Mooney-Rivlin \eqref{eq:mooney}, Gent \eqref{eq:gent} and Haine-Wilson \eqref{eq:HW}, respectively, for the quantity $Q$. The model error corresponds to the relative error between a very fine FEM solution $\bu_{\mbox{fine}}$ and the experimental quantity of interest. The discretization error corresponds to the relative error between the computed solution on the current mesh and the computed solution on a very fine mesh, i.e.
$$
\mbox{model error}=\dfrac{|Q(\bu_{\mbox{fine}})-20|}{20}
\mbox{ and discr. error}=\dfrac{|Q(\bu_h)-Q(\bu_{\mbox{fine}})|}{Q(\bu_{\mbox{fine}})}.
$$
The parameters of each law have been estimated from \cite{meunier2008mechanical} and are provided in Table \ref{tab:parameter_law}.
We give in Figure \ref{fig:deformed mesh} the refined mesh in the case of the Haine-Wilson law (left) and the deformed geometry when we apply the load (right). We remark that the refinement occurs mainly on the top the silicone where is localised the quantity of interest but also near the holes.

 \begin{table}[]
        \centering
\begin{tabular}{llll}\hline
Mooney&$C_{10}=0.14$&$C_{01}=0.023$&\\
Gent&$E=0.97$&$J_m=13$&\\
Haines-Wilson&$C_{10}=0.14$&$C_{20}=-0.0026$&$C_{30}=0.0038$\\
&$C_{01}=0.033$&$C_{02}=0.00095$&$C_{11}=-0.0049$\\\hline
        \end{tabular}       
        \caption{First test case (silicone sample). Values of the constitutive parameters  of each hyperelastic law, following \cite{meunier2008mechanical}.}
        \label{tab:parameter_law}
    \end{table}

    \begin{table}[]
        \centering
\begin{tabular}{|c|c|c|c|c|c|c|c|c|}\hline
\multicolumn{3}{|c}{Mooney-Rivlin}&\multicolumn{3}{|c|}{Gent}&\multicolumn{3}{c|}{Haine Wilson}\\
\hline
   Nb& Model&Discr.&Nb& Model&Discr.&Nb& Model&Discr. \\
     Cells& Error&  Error&Cells& Error& Error&Cells& Error& Error\\\hline
1333&&7.0\%&1333&&9.6\%&  1333   &&8.0\% \\
1664&&8.1\%&1654&&10.4\%&  1706  &  &8.8\%\\
3047&&7.2\%&3281&&9.1\%&  3441 &  &7.8\% \\
5307&17.0\%&4.3\%&5414&2.1\%&5.9\%&  5692 &0.5\%&4.6  \%\\
9426&&3.0\%&9947&&3.7\%&   10383&&2.9 \%\\
18118&&2.3\%&18088&&2.6\%&   19489&&2.0\%\\
33982&&1.9\%&34883&&2.0\%&   38466 &&1.8\%\\\hline
        \end{tabular}       
        \caption{First test case (silicone sample).  Relative model and discretisation error of the quantity of interest $Q_1$ with respect to the number of cells. The model error is computed thanks to the experimental data. 
        }
        \label{tab:law}
    \end{table}

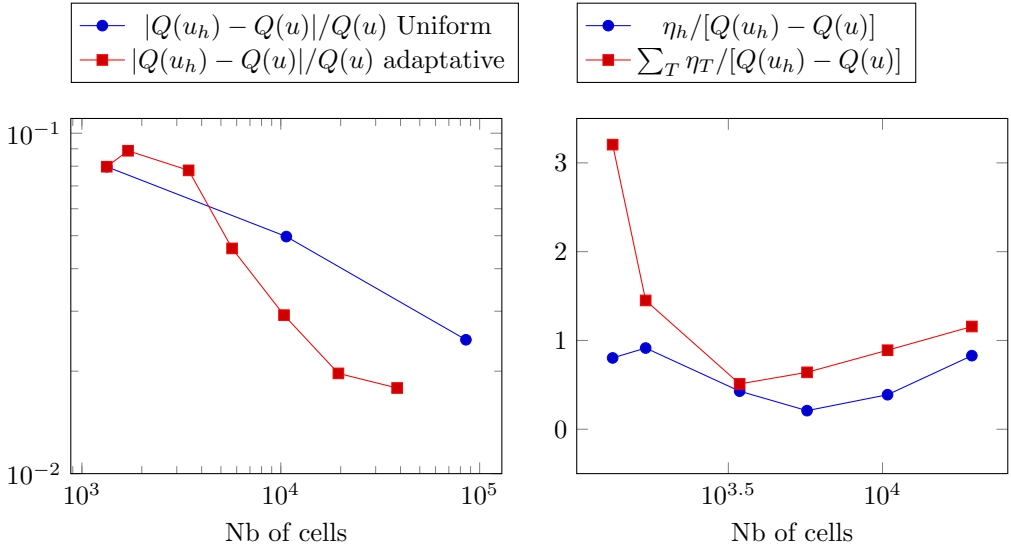
\begin{figure}
\centering
\begin{tikzpicture}\begin{loglogaxis}[
width = .49\textwidth, xlabel = Nb of cells, ymin=1e-2,
            legend style = { at={(0,1.1)},anchor=south west, legend columns =1,
			/tikz/column 2/.style={column sep = 10pt}}]

\addplot coordinates { 
(1333.0,0.07974384496146314)
(10664.0,0.0497004842772945)
(85312.0,0.024731635079811583)
};
\addplot coordinates{
(1333.0,0.07974384496146314)
(1706.0,0.08881530907994144)
(3441.0,0.07780126371075838)
(5692.0,0.045881086154962114)
(10383.0,0.0292177412714193)
(19486.0,0.019705163641890407)
(38466.0,0.01784846789296601)
};
\legend{$|Q(u_h)-Q(u)|/Q(u)$ Uniform,$|Q(u_h)-Q(u)|/Q(u)$ adaptative}
\end{loglogaxis}
\end{tikzpicture}
\quad
\begin{tikzpicture}\begin{semilogxaxis}[width = .49\textwidth, xlabel = Nb of cells, ymin=-0.5,ymax=3.5,
            legend style = { at={(0,1.1)},anchor=south west, legend columns =1,
			/tikz/column 2/.style={column sep = 10pt}}]

\addplot coordinates { 
(1333.0,0.8030303154863103)
(1706.0,0.9144736108669475)
(3441.0,0.4287095121982907)
(5692.0,0.2092797043313161)
(10383.0,0.38865031573267556)
(19486.0,0.8287531974788899)
};
\addplot coordinates{
(1333.0,3.2053004788007224)
(1706.0,1.4503989000584685)
(3441.0,0.510507948303809)
(5692.0,0.6409859330065689)
(10383.0,0.8898992292127493)
(19486.0,1.1567722834967615)
};
\legend{$\eta_h/[Q(u_h)-Q(u)]$,$\sum_T\eta_T/[Q(u_h)-Q(u)]$}
\end{semilogxaxis}
\end{tikzpicture}
\caption{First test case (silicone sample). Relative error of discretisation (left) and efficiency of the estimator (right).}
\label{fig:errorH1_manufactured}
\end{figure}


\subsection{Second test case : Human heel} 



A pressure ulcer (PU) is a wound stemming from excessive loads on biological soft tissues which leads to ischemia, which in turn triggers tissue necrosis. Two fifth of the patients taken in charge by a reanimation unit or in a geriatric unit will develop a PU. 40\% of those ulcers are located at the posterior part of the heel because the patient stays for hours lying on his back without moving \cite{Perneger1998}. This condition is often followed by amputation of part of the foot. The high prevalence of the pathology is the motivation behind the recent development of PU prevention strategies. Some of these approaches resort to personalized biomechanical modeling of the patient's soft tissues, where tissue compression is numerically predicted based on the loads measured underneath the bedridden patient's heels. The current consensus in the PU prevention community is that PU risk assessment should be based on an indicator of tissue suffering derived from the Von Mises stress, see \cite{loerakker2011effects}. Thus, a personalized biomechanical model should be able to predict the onset of a PU by continuously monitoring this quantity of interest, and by triggering an alarm if the risk exceeds a pre-defined threshold.


The accuracy of such prediction not only highly depends on the accuracy of the determined mechanical properties of the heel tissue, on the correct boundary conditions used in the simulation, but also on the numerical method (here the FEM) which is capable to solve the problem in a way that the error is controlled.

\begin{figure}[!htbp]
\centering
\includegraphics[width=0.4\columnwidth]{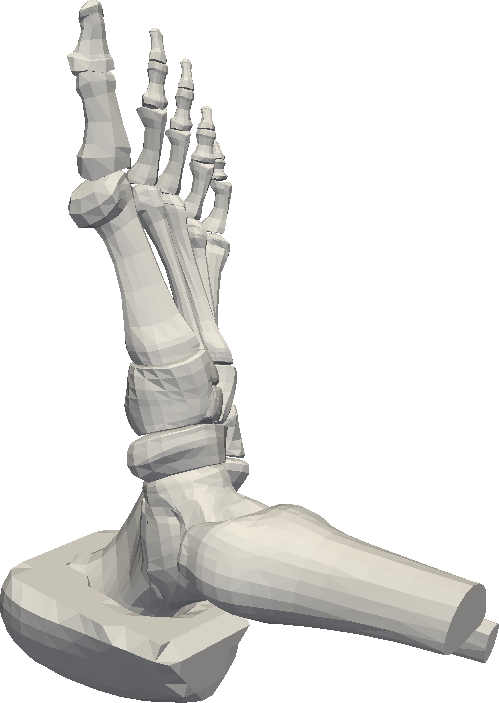}
\caption{A part of a heel tissue model used in simulations in which its orientation corresponds to the situation when the patient is on bed.}
\label{fig:nonlinear/heel_model}
\end{figure}

\begin{figure}[!htbp]
 \centering
       \begin{subfigure}[b]{0.4\textwidth}
	  \centering
      \includegraphics[width=1\columnwidth]{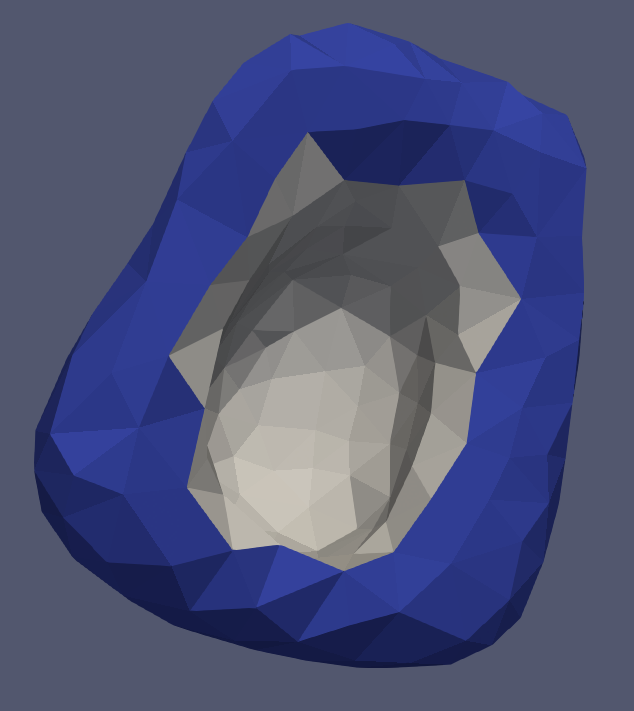}
	  \caption{}
	  \label{fig:heel_fixed_surface}
      \end{subfigure}%
            ~ 
      \begin{subfigure}[b]{0.4\textwidth}
	  \centering
	  \includegraphics[width=1\columnwidth]{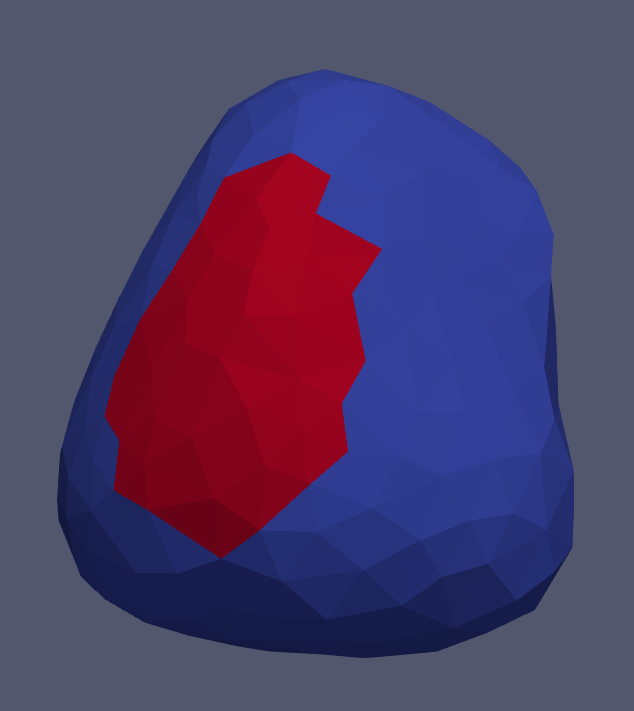}
	  \caption{}
	  \label{fig:heel_neumann_and_region_interest}
      \end{subfigure}%
\caption{The heel tissue is considered to be fixed on the surface which has contact with the calcaneum, and on the upper surface, shown by gray colour in (\subref{fig:heel_fixed_surface}); the tissue surface where a pressure is applied is shown by red colour, whereas a region of interest is also shown by red colour (\subref{fig:heel_neumann_and_region_interest}).}\label{fig:heel_boundary_conditions}
\end{figure}

\cref{fig:nonlinear/heel_model} shows a heel tissue model used in our simulation. Its orientation corresponds to the situation when the patient is on bed. Boundary conditions used are shown in \cref{fig:heel_boundary_conditions}. We consider Von Mises stress is a good factor to predict the tissue damage. A region in which we consider the heel tissue is vulnerable is shown by cyan colour in \cref{fig:heel_neumann_and_region_interest}. Our quantity of interest is thus expressed through the first Piola-Kirchhoff stress tensor $\bPi$ over the domain of interest $\omega$, as
\begin{equation}
Q(\bPi(\bu,\bpp)) = \int_\omega \sqrt{\frac{1}{2} \left( (\bPi_{11}-\bPi_{22})^2 + (\bPi_{22}-\bPi_{33})^2 + (\bPi_{33}-\bPi_{11})^2 +3(\bPi_{12}^2 +\bPi_{21}^2 + \bPi_{23}^2+ \bPi_{32}^2 + \bPi_{31}^2+ \bPi_{13}^2)\right)}.
\end{equation}
%
The heel tissue is supposed to behave like an 
incompressible hyperelastic material with the mechanical properties: 
Mooney-Rivlin with $C_{10}=16.6 kPa$ and $C_{01}=0$ (see \cite{perrier:hal-01930207}).

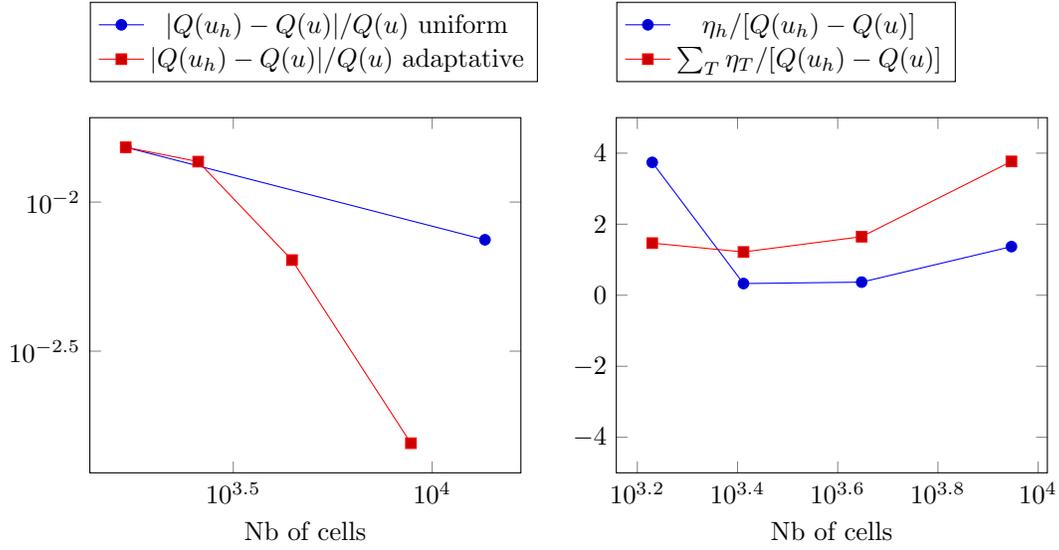
\begin{figure}
\centering
\begin{tikzpicture}\begin{loglogaxis}[
width = .49\textwidth, xlabel = Nb of cells,
            legend style = { at={(0,1.1)},anchor=south west, legend columns =1,
			/tikz/column 2/.style={column sep = 10pt}}]

\addplot coordinates { 
(1697.0,0.015271988367355122)
(13576,0.0074748147418371704)
};
\addplot coordinates{
(1697.0,0.015271988367355122)
(2581.0,0.01366912947814863)
(4443.0,0.006390160766408987)
(8843.0,0.0015532929202563607)
};

\legend{$|Q(u_h)-Q(u)|/Q(u)$ uniform,$| Q(u_h)-Q(u)|/Q(u)$ adaptative}
\end{loglogaxis}
\end{tikzpicture}
\quad
\begin{tikzpicture}\begin{semilogxaxis}[width = .49\textwidth, xlabel = Nb of cells, ymin=-5,ymax=5,
            legend style = { at={(0,1.1)},anchor=south west, legend columns =1,
			/tikz/column 2/.style={column sep = 10pt}}]

\addplot coordinates { 
(1697.0,3.7404742866423697)
(2581.0,0.3288313980969957)
(4443.0,0.36866063970302765)
(8843.0,1.3687666355098917)
};
\addplot coordinates{
(1697.0,1.4678614343464433)
(2581.0,1.2179865385206503)
(4443.0,1.648423312436307)
(8843.0,3.7673339801659944)
};
\legend{$\eta_h/[Q(u_h)-Q(u)]$,$\sum_T\eta_T/[Q(u_h)-Q(u)]$}
\end{semilogxaxis}
\end{tikzpicture}
\caption{Second test case (Human heel undergone surface pressure). Relative error of discretisation (left) and efficiency of the estimator (right).}
\label{fig:nonlinear/heel_error_norms}
\end{figure}

\begin{figure}[!htbp]
 \centering
       \begin{subfigure}[b]{0.25\textwidth}
	  \centering
      \includegraphics[width=1\columnwidth]{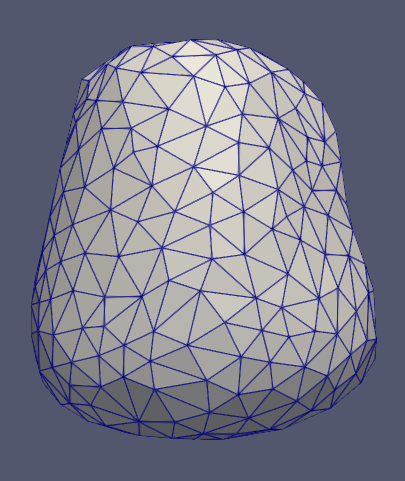}
	  \caption{}
	  \label{fig:heel_mesh_back_iter0}
      \end{subfigure}%
            ~ 
      \begin{subfigure}[b]{0.25\textwidth}
	  \centering
	  \includegraphics[width=1\columnwidth]{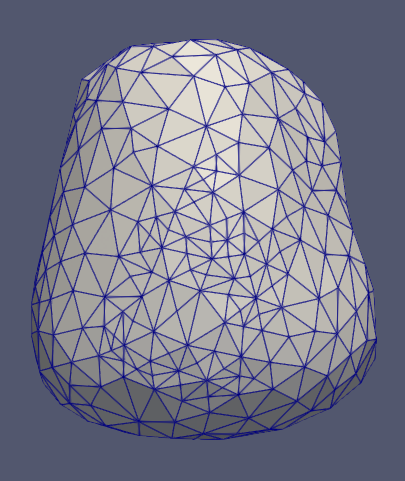}
	  \caption{}
	  \label{fig:heel_mesh_back_iter1}
      \end{subfigure}%
       ~ 
      \begin{subfigure}[b]{0.25\textwidth}
	  \centering
	  \includegraphics[width=1\columnwidth]{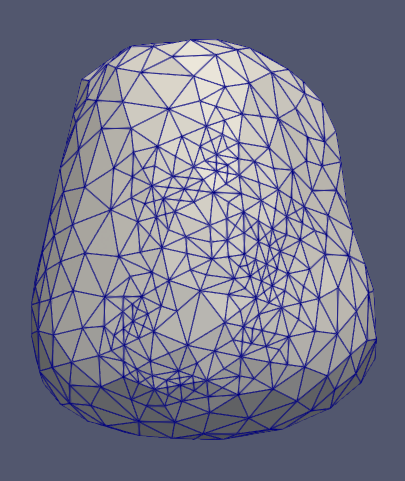}
	  \caption{}
	  \label{fig:heel_mesh_back_iter2}
      \end{subfigure}%
             ~ 
      \begin{subfigure}[b]{0.25\textwidth}
	  \centering
	  \includegraphics[width=1\columnwidth]{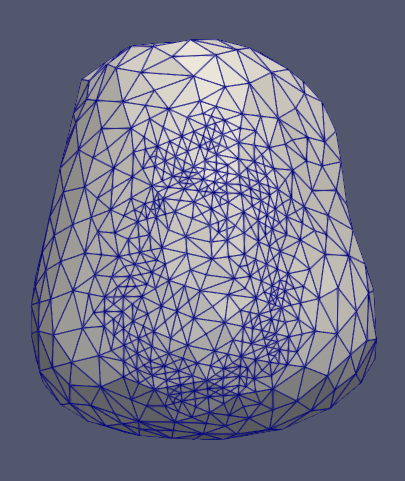}
	  \caption{}
	  \label{fig:heel_mesh_back_iter3}
      \end{subfigure}%
      \\
             \begin{subfigure}[b]{0.25\textwidth}
	  \centering
      \includegraphics[width=1\columnwidth]{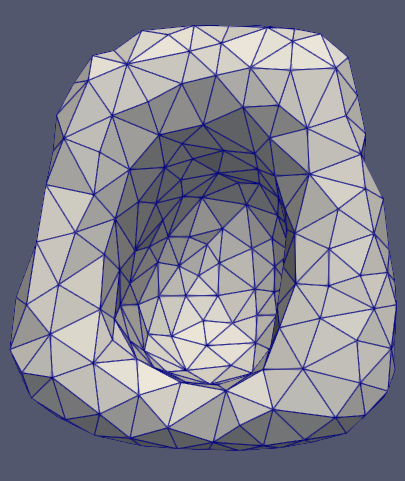}
	  \caption{}
	  \label{fig:heel_mesh_front_iter0}
      \end{subfigure}%
            ~ 
      \begin{subfigure}[b]{0.25\textwidth}
	  \centering
	  \includegraphics[width=1\columnwidth]{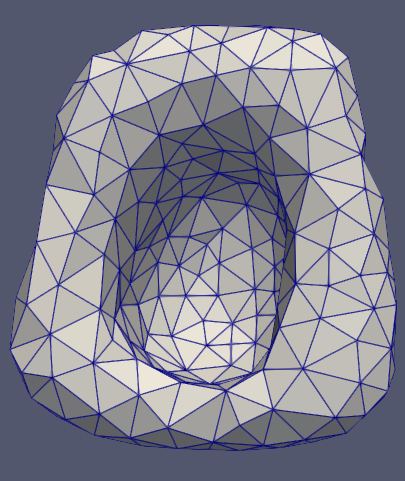}
	  \caption{}
	  \label{fig:heel_mesh_front_iter1}
      \end{subfigure}%
       ~ 
      \begin{subfigure}[b]{0.25\textwidth}
	  \centering
	  \includegraphics[width=1\columnwidth]{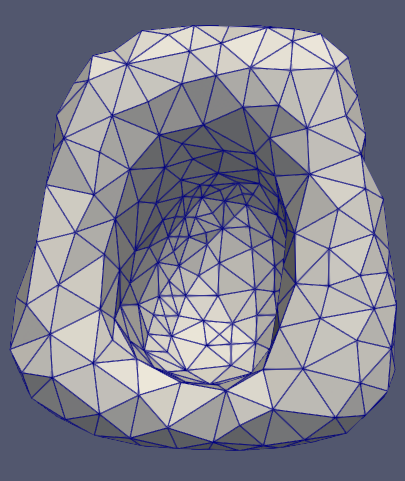}
	  \caption{}
	  \label{fig:heel_mesh_front_iter2}
      \end{subfigure}%
             ~ 
      \begin{subfigure}[b]{0.25\textwidth}
	  \centering
	  \includegraphics[width=1\columnwidth]{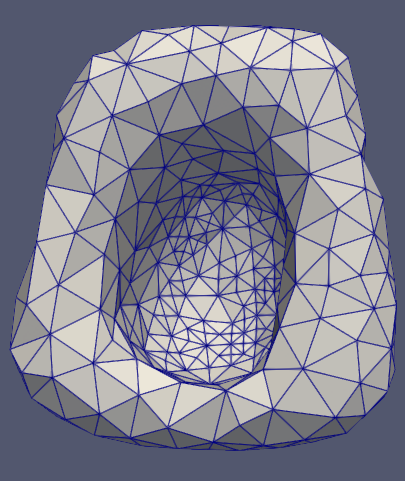}
	  \caption{}
	  \label{fig:heel_mesh_front_iter3}
      \end{subfigure}%
\caption{Second test case (Human heel undergone surface pressure). From the back, initial mesh (\subref{fig:heel_mesh_back_iter0}), adaptive mesh obtained after the first iteration of refinement (\subref{fig:heel_mesh_back_iter1}), after the second iteration of refinement (\subref{fig:heel_mesh_back_iter2}), the final mesh obtained after the $3th$ iteration of refinement (\subref{fig:heel_mesh_back_iter3}). The same meshes from the front (\subref{fig:heel_mesh_front_iter0})-(\subref{fig:heel_mesh_front_iter3}).}\label{fig:heel_adaptive_meshes}
\end{figure}

\cref{fig:nonlinear/heel_error_norms} shows the convergence of the error estimator under both uniform and adaptive refinements. It is observed again that using adaptive refinement scheme is more advantageous than the uniform one since the corresponding error converges with higher rate. Refinement patterns are shown in \cref{fig:heel_adaptive_meshes}.
 

We remark that the algorithm refine in the lower part of the heel which is in coherence with \cite{van2020risk,GEFEN2010124,SHAULIAN2021104261,luboz2015influence}.
These zones correspond to the onset of pressure ulcers (see \cite{traa2018mri}).
The added value of our DWR-driven adaptative meshing is that it automatically refines according to a quantity of interest. 
It not only optimises locally in the region of interest but also takes into account the far-field / global errors that contribute. 

\section{Discussion}
\label{sec:Conclusions}

We summarize below our main achievements, discuss their current limitations and suggest some perspectives.

\subsection{Main achievements}

We illustrate in this paper the feasibility of carrying out mesh refinement in modern finite element environments such as FEniCS, driven by a goal oriented error estimate, to improve the accuracy of an arbitrary user-defined quantity of interest. This is done for a three-dimensional nonlinear problem that involves incompressible hyperelasticity. This model is of common use in the biomechanics community. 
Notably, large displacements and large strain are taken into account. Furthermore, there is no difficulty to incorporate incompressibility. 

The prediction outcome of the simulations has been confronted to experimental measurements on silicone sheets. Validation and verification of simulations in computational biomechanics is a complicated issue, since measurements on real patients and real clinical situations lack of precision and reproductibility. Indeed, some difficulties arise to control and know accurately all the relevant parameters. Though they are of course not identical as human soft tissues, materials such as silicone are interesting, since the numerical model can be calibrated with precision and numerical simulations can be confronted to measurements. 
However, it is still very difficult to find some published material for this purpose, and we took advantage of the study of Meunier et al \cite{meunier2008mechanical}.
In addition, 
let us emphasize the following points


\begin{enumerate}
    \item The test with the silicone sample provides extra information about the modelling error, in the sense it also quantifies the predictive power of the constitutive law.
    
    \item When the mesh is too coarse, the discretization error is of the same magnitude of this modelling error. 
    
\item With adaptive mesh refinement, the discretization error can be controlled and driven below a given threshold that makes it negligible, without having to overrefine the mesh.

    \item Another example close to clinical biomechanics has been carried out on a three-dimensional complex patient-specific geometry of the heel. 
    
    \item The method is easy to implement in an environment such as FEniCS, and the scripts are freely available. They can be transposed without much difficulty to other similar environments such as SONICS\cite{mazier2022sonics}, GetFEM++\cite{renard2020getfem}, FreeFEM++\cite{hecht2012new} or SciKit-FEM\cite{gustafsson2020scikit}.

    \item Notably the automatic differentiation tools now available in modern finite element software facilitate a lot the assembly of the dual problem. Solving the, linear, dual problem remains inexpensive in comparison with the total solution procedure needed for the nonlinear problem. Last but not least, the solution of this dual problem is the basis of counterintuitive refinement strategies, much more efficient than adhoc refinement.

\end{enumerate}

\subsection{Current limitations}

Let us point out as well some limitations of the proposed methodogy:

\begin{enumerate}
    \item The dual solution needs to be approximated, and when the same finite element spaces are used for the primal and dual problem, an extra extrapolation step needs to be carried out. The accuracy of this procedure may be quantified more precisely and improved. In practice it has revealed small effect and does not hamper the efficiency of the methodology, but better results may even be expected if this point is improved. Another possibility is to use higher order spaces for the dual problem, but this solution is much more expensive.
    
    \item The adjoint problem, though it is linear, inherits its coefficients from the nonlinear primal problem. In specific situations, it may be ill conditionned and this issue would need further investigation.
    
    
    \item From the representation formula \eqref{eq:representation_DWR} of Becker and Rannacher at the core of the error estimate, we neglect the linearization error. Of course, this one is complicated to estimate in general and is expected to be small, but this point may deserve to be studied more carefully in the future.

    \item Reference solutions are computed solutions on fine meshes or eventually experimental data. It could be interesting also to test the methodology with manufactured solutions in hyperelasticity \cite{payan-ohayon-2017,blaise:2022,chamberland2010comparison}.
    
    \end{enumerate}

\subsection{Perspectives}

To decrease the computation time, it can be interesting to perfom the refinement at each step throughout the loading. It can be significant when considering non-linear problems.
A stimulating perspective would be first to combine the current methodology with techniques for model selection, and to estimate more systematically the model error. Also, since, for patient-specific biomechanics, some data can be used for (possibly) real time parameter calibration, it would be interesting to take advantage of the flexibility of the current framework to account for parameter calibration, as already done in Becker \& Vexler \cite{becker2004,becker2005} for a general setting. Another point would consist in making the methodology available in software that are of common use in the whole biomechanics community. 
More generally, a perspective can be the comparison of the model et discretisation error with the errors coming from the geometry, the parameters of the model and the forces applied on the organ.





\section{Acknowledgements}

The authors thank Roland Becker, Jacques Ohayon, Yohan Payan and Yves Renard for their advices. They also thank the AMIES for its support, particularly Magalie Fredoc, Antoine Lejay and Christophe Prudhomme.

\newpage
\appendix

\section{Expression of the estimator and algorithms}

We give below for each cell-wise contribution:
\begin{equation}
\eta_K = \left| \int_K \bR_{\bu} \cdot (\bz_h-I_h(E_h(\bz_h)) \mathrm{d}\Omega + \int_K \bR_{\bpp} \cdot (\bpp_h-I_h(E_h(\bpp_h)) \mathrm{d}\Omega +  \int_{\partial K} \bJ \cdot (\bz_h-I_h(E_h(\bz_h)) \mathrm{d}\gamma   \right|
\end{equation}
with, $E_h$ and $I_h$ are resp. the extrapolation (see \cite{rognes-logg-2013}) and the Lagrange interpolation, the interior residual
\begin{equation*}
\bR_{\bu} = \bB + 
\mathrm{div} \, \bPi(\bu_h)
\mbox{ and }\bR_{\bpp}=\det(\boldsymbol{C})-1
\end{equation*}
and the stress jump 
\begin{equation*}
\bJ = \begin{cases}
   - \frac{1}{2}   [[\bPi(\bu_h)]]
       & \text{ if } \gamma \not \subset \Gamma, \\
      \bT -  \bPi(\u_h)  \cdot \bN                                 & \text{ if } \gamma \subset \Gamma_N, \\
      \bzero                                                                     & \text{ if } \gamma \subset \Gamma_D.  
      \end{cases}
\end{equation*}
The jump can be defined for a function $\bv_h$ on a facet $F$ between two cells $K$ and $K'$ by 
$[[\bv_h]]=\bv_{h|K}\cdot n_K + \bv_{h|K'}\cdot n_{K'}$, where $n_K$ and $n_{K'}$ are the normal of $K$ and $K'$ on $F$.
\begin{algorithm}
\caption{Algorithm for mesh refinement}\label{algo:adaptivity}
\begin{algorithmic}
\State Select an initial triangulation $mesh_i$ of the domain $\Omega$
\State Guest solution $(\bu_i^{(0)},\bpp_i^{(0)})$
\While{$ \sum_K \eta_K > \epsilon$}
\State $F(\bu_i,\bpp_i;\bv,\bq) \gets A(\bu_i,\bpp_i;\bv,\bq) - L(\bv,\bq)$
\State $\bu_i,\bpp_i \gets$ NewtonMethod $(F(\bu_i,\bpp_i;\bv,\bq),(\bu_i^{(0)},\bpp_i^{(0)}))$
\Comment{Problem \eqref{eq:primal_weak}, see Algo \ref{algo:newton_method}} 
\State $\bz_i ,\bw_i\gets$ DualProblem $(\bu_i,\bpp_i, Q)$                      \Comment{Problem \eqref{eq:discrete_dual_problem}, see Algo \ref{algo:dual_problem}} 
\State $\eta_K \gets$ ComputeEstimator $(\bu_i,\bpp_i,\bz_i,\bw_i)$
\State markedElements $\gets$ DorflerMarking $(\eta_K, \alpha)$ \Comment{See Algo \cref{algo:Dorfler_making}}
\State $mesh_i \gets mesh_i$.refine(markedElements)
\State Compute $\sum_K \eta_K$
\EndWhile
\end{algorithmic}
\end{algorithm}
%
\begin{algorithm}
\caption{Solving a non-linear problem: NewtonMethod $(F(\bu_i^{(0)},\bv),\bu_i^{(0)})$}\label{algo:newton_method}
\begin{algorithmic}
\State $(\bu_k,\bpp_k) = (\bu_i^{(0)},\bpp_i^{(0)}) $
\While{$\lvert (\delta \bu,\delta\bpp) \rvert > \epsilon$}
 \State $ F'(\bu_k ; \delta \bu , \bv)  = -F(\bu_k; \bv) $    \Comment{Solve for $\delta \bu$} 
 \State $ (\bu_{k+1},\bpp_{k+1}) \gets (\bu_k + \delta \bu,\bpp_k + \delta \bpp) $    \Comment{Update the solution}
 \State Compute $\lvert (\delta \bu,\delta\bpp) \rvert$
\EndWhile
\end{algorithmic}
\end{algorithm}
\begin{algorithm}
\caption{Solving the dual problem: DualProblem $(\bu_i,Q)$}\label{algo:dual_problem}
\begin{algorithmic}
\State Compute $A'(\bu_i,\bpp_i| \bz_i,\bw_i;\bv,\bq) $
\State Compute $Q'(\bu_i,\bpp_i; \bv,\bq)$
\State $(\bz_i,\bw_i) \gets$ solve $(A'(\bu_i,\bpp_i| \bz_i,\bw;\bv,\bq) = Q'(\bu_i,\bpp_i; \bv,\bq) )$ \Comment{Solve the linear system}
\end{algorithmic}
\end{algorithm}
\begin{algorithm}
\caption{Mark elements after D\"orfler strategy by providing a element-wise estimator $\eta_K = [\eta_{K_1}, \eta_{K_2}, \dots \eta_{K_N}]$, and $0 < \alpha < 1$ a parameter which characterises the marking rate: the smaller the value of $\alpha$ is, the fewer the number of elements will be marked: DorflerMarking $(\eta_K, \alpha)$ }\label{algo:Dorfler_making}
\begin{algorithmic}
\State Sort the elements $K_i$ after descending order of the corresponding estimator $\eta_{K_i}$
\State Mark the first $M$ elements such that 
\[\mathrm{markedElements} \gets \min \left \{ M \in N \:\middle |\: \sum\limits_{i=1}^{M}\eta_{K_i}\geq\alpha\sum_{i=1}^N \eta_{K_i} \right \}. \]
\end{algorithmic}
\end{algorithm}

\bibliographystyle{abbrv}
\bibliography{MainManuscript.bib}

\end{document}